\documentclass{amsart}\usepackage[reset,a4paper,vmargin=5truecm,hmargin=4truecm]{geometry}
\usepackage{amsmath,amsthm,amssymb}

\usepackage{graphicx, color}
\usepackage{url}
\usepackage{cite}
\usepackage{hyperref}

\newtheorem{thm}{Theorem}

\begin{document}

\title[Convex Polygons for Aperiodic Tiling]{Convex Polygons for Aperiodic Tiling}

\address{The Interdisciplinary Institute of Science, Technology and Art\\
Suzukidaini-building 211, 2-5-28 Kitahara, Asaka-shi, Saitama, 351-0036, 
Japan}
\email{ismsugi@gmail.com}
\author{Teruhisa Sugimoto}
\maketitle	

\begin{abstract}
If all tiles in a tiling are congruent, the tiling is 
called monohedral. Tiling by convex polygons is called edge-to-edge if any 
two convex polygons are either disjoint or share one vertex or one entire 
edge in common. In this paper, we prove that a convex polygon that can generate an 
edge-to-edge monohedral tiling must be able to generate a periodic tiling.
\end{abstract}

\section{Introduction}
\label{section1}

A \textit{tiling} of the plane is an exact covering of the plane by a collection of 
sets without gaps or overlaps (except for the boundaries of the sets). More 
precisely, a collection of sets (the ``tiles'') is a tiling of the plane if 
their union is the entire plane, but the interiors of different tiles are 
disjoint. The tiles are frequently polygons and are all congruent to one 
another or at least congruent to one of a small number of 
\textit{prototiles}~\cite{G_and_S_1987, Hallard_1991}. 
If all tiles in the tiling are of the same size and shape (i.e., congruent), the 
tiling is called \textit{monohedral }~\cite{G_and_S_1987}. Therefore, a 
prototile of a monohedral tiling by polygons is single polygon in the tiling, 
and we call the polygon the \textit{polygonal tile}~\cite{G_and_S_1987, 
Sugimoto_2012a, Sugimoto_2012b, Sugimoto_2015, Sugimoto_III, 
Sugimoto_NoteTP, Sugi_Ogawa_2006}.

A tiling of the plane is \textit{periodic} if the tiling can be translated onto itself 
in two nonparallel directions. More precisely, a tiling is periodic if it coincides 
with its translation by a nonzero vector. A set of prototiles is called 
\textit{aperiodic} if congruent copies of the prototiles admit infinitely many 
tilings of the plane, none of which are periodic. It must be emphasized 
that no periodic tilings are permitted at all, even using just one of the 
prototiles~\cite{Akiyama_2012, G_and_S_1987, Hallard_1991, Socolar_2011}. 
(A tiling that has no periodicity is called nonperiodic. On the other 
hand, a tiling by aperiodic sets of prototiles is called aperiodic. Note 
that, although an aperiodic tiling is a nonperiodic tiling, a nonperiodic 
tiling is not necessarily an aperiodic tiling.) The Penrose tiling, which is 
known as a quasiperiodic tiling, is a nonperiodic tiling, and it can also be 
considered an aperiodic tiling that is generated by the aperiodic set of 
prototiles with a matching condition (see Section~\ref{section2} for details).

For aperiodic tilings, the critical problem is to find sets of aperiodic 
prototiles that are essentially different from the ones already known. In 
particular, there are following two problems~\cite{Hallard_1991}.

\bigskip
\noindent
(i) Is there a single aperiodic prototile (with or without a matching 
condition), that is, one that admits only aperiodic tilings by congruent 
copies?

\noindent
(ii)  It is well known that there is a set of three convex polygons that are 
aperiodic with no matching condition on the edges. Is there a set of 
prototiles with a size less than or equal to two that is aperiodic?

\bigskip
\noindent
Problem (i) is called the Einstein Problem (ein stein = one stone), and the 
affirmative solution (a convex hexagon with a matching condition) was shown in 
recent years~\cite{Socolar_2011}. The elements of an aperiodic set of three convex 
polygons with no matching condition on the edges in Problem (ii) are one 
convex hexagon and two convex pentagons~\cite{G_and_S_1987, 
Hallard_1991}. Then, the aperiodic tiling by the three convex polygons 
is based on the Penrose tiling (see Subsection~\ref{subsection2_2}).

Tiling by convex polygons is called \textit{edge-to-edge} if any two convex 
polygons are either disjoint or share one vertex or one entire edge in 
common~\cite{G_and_S_1987, Sugimoto_2012a, Sugimoto_2012b, 
Sugimoto_2015, Sugimoto_III, Sugimoto_NoteTP, Sugi_Ogawa_2006}. 
We have studied the convex pentagonal tiles that can generate an 
edge-to-edge tiling. As the result, we find the following~\cite{Sugimoto_2012b, 
Sugimoto_NoteTP}.

\begin{thm}\label{t2}
Without matching conditions other than ``edge-to-edge,'' no 
single convex polygon can be an aperiodic prototile.
\end{thm}

In Section~\ref{section2}, we explain the aperiodic set of prototiles from the 
Penrose tiling and introduce the aperiodic tiling by one convex hexagon and two 
convex pentagons with no matching condition on the edges. In Section~\ref{section3}, 
we present the proof of Theorem~\ref{t2}. Section~\ref{section4} presents an 
estimate for the aperiodic tiling by a convex polygonal tile from the known facts 
on convex polygonal tiles.

\section{Aperiodic Tiling and Aperiodic Set of Prototiles}
\label{section2}

\subsection{Penrose Tiling}
\label{subsection2_1}

Here, an aperiodic set of prototiles is introduced using the Penrose tiling. 
This is the topic relevant to Problem (i) on aperiodic tiling mentioned in 
Section~\ref{section1}.

\begin{figure}[htbp]
 \centering\includegraphics[width=12.8cm,clip]{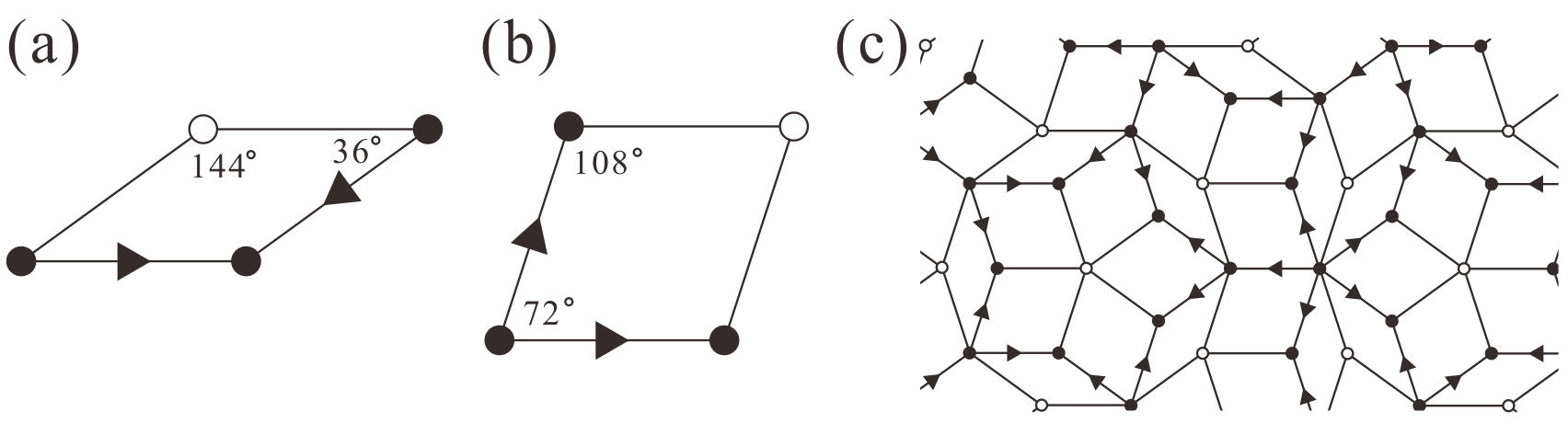} 
  \caption{{\small 
Penrose tiling and Penrose rhombuses of an aperiodic set.} 
\label{fig1}
}
\end{figure}

\begin{figure}[htbp]
 \centering\includegraphics[width=12.8cm,clip]{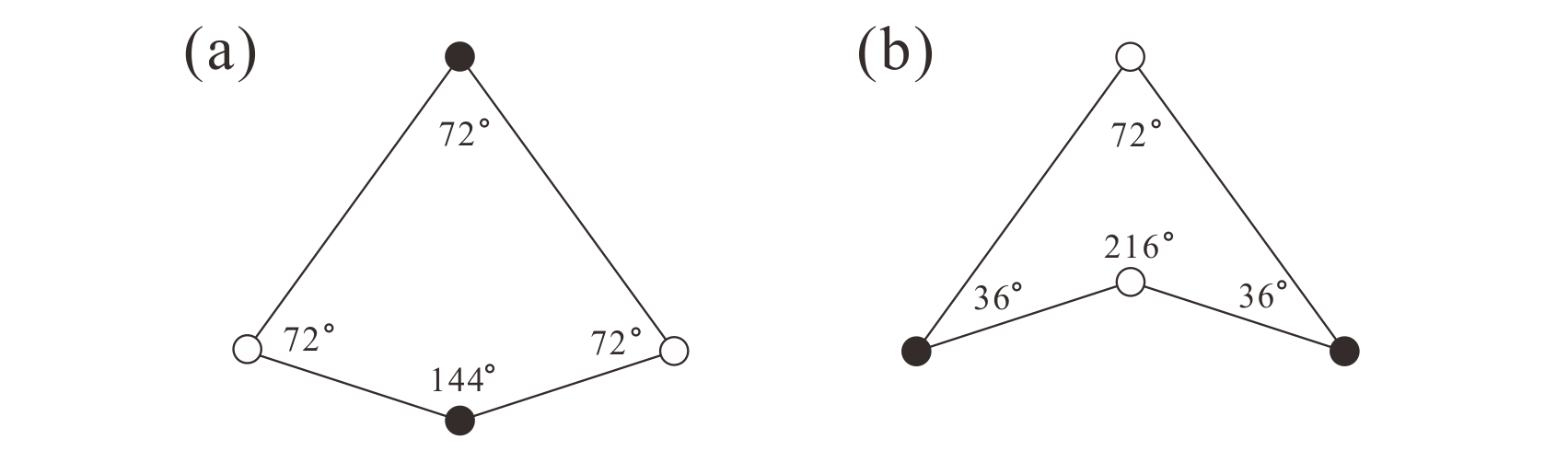} 
  \caption{{\small 
The kite and the dart of the Penrose aperiodic set.} 
\label{fig2}
}
\end{figure}

The rhombuses in Figure~\ref{fig1}(a) and (b) (called Penrose rhombuses) are the 
prototiles of the Penrose tiling. We note that all vertices of the two 
rhombuses are colored; there exist edges with orientations, and the length 
of the edges of the rhombus in Figure~\ref{fig1}(a) is equal to that of the rhombus 
in Figure~\ref{fig1}(b). To obtain the tiling that is generated by the rhombuses in 
Figure~\ref{fig1}(a) and (b), the vertices always meet with the same color; it is 
edge-to-edge, and the edges with orientations must match the direction of 
the orientations. The way of matching that is required for the generation of 
the tiling is called a \textit{matching condition}~\cite{G_and_S_1987, 
Hallard_1991}. The set of the two rhombuses in Figure~\ref{fig1}(a) and (b) is an 
aperiodic set of prototiles with a matching condition ~\cite{G_and_S_1987}. 
Therefore, the Penrose tiling in Figure~\ref{fig1}(c) is an aperiodic tiling in which 
the rhombuses in Figure~\ref{fig1}(a) and (b) are generated according to the 
matching condition.  In addition, there is another way of choosing the two 
prototiles of the Penrose tiling: the kite (see Figure~\ref{fig2}(a)) and the dart (see 
Figure~\ref{fig2}(b)). The edges of the kite and the dart have two lengths in the ratio 
$1:\tau$, where ${\tau = \left( {1 + \sqrt 5 } \right)} \mathord{\left/ 
{\vphantom {{\tau = \left( {1 + \sqrt 5 } \right)} 2}} \right. 
\kern-\nulldelimiterspace} 2 \approx 1.1618\ldots $ is the golden number. 
The vertices are colored with two colors, say black and white, as shown. To 
obtain the Penrose tiling by using the kite and the dart, we must place equal 
edges together and also match the colors at the vertices. The set of the 
kite and the dart in Figure~\ref{fig2} is an aperiodic set of prototiles with a matching 
condition ~\cite{G_and_S_1987, Hallard_1991}. Note that each prototile in 
Figures~\ref{fig1} and \ref{fig2} cannot independently generate a tiling when 
it is placed according to a matching condition.

\renewcommand{\figurename}{{\small Figure.}}

\renewcommand{\figurename}{{\small Figure.}}

In the above explanation for a matching condition, we use vertices with 
color conditions. If the edges are admitted to use a color, the vertices 
do not need to have a matching condition. 
That is, a matching condition is a condition concerning edge 
matching. It must be emphasized that the matching condition can be 
represented by assigning colors and orientations to some of the 
edges of the prototile.

\begin{figure}[htbp]
 \centering\includegraphics[width=12.8cm,clip]{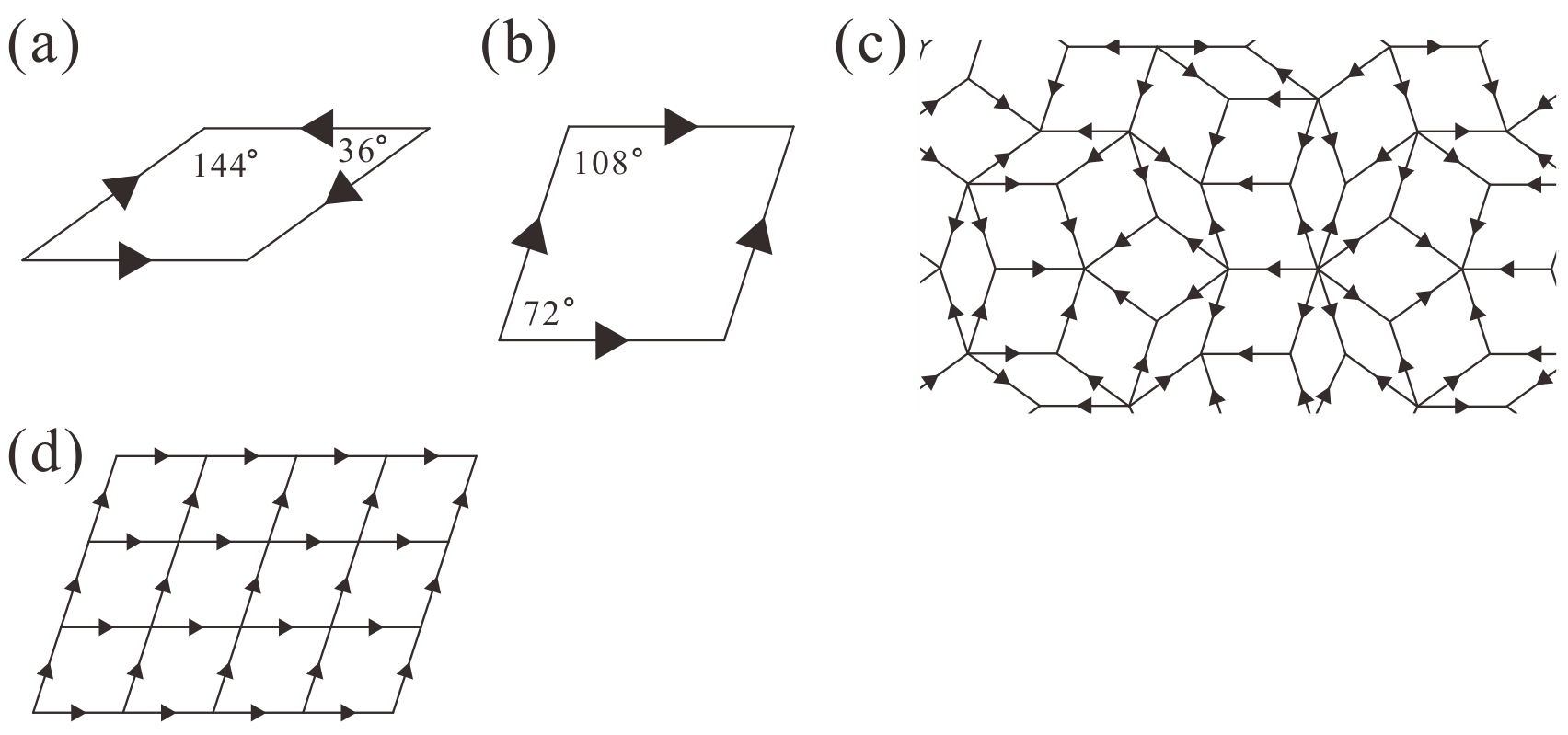} 
  \caption{{\small 
Prototiles of the Penrose tiling, which is not an aperiodic set.}
\label{fig3}
}
\end{figure}

Next, in order to understand the aperiodic set of prototiles, the prototiles 
of the Penrose tilings shown in Figure~\ref{fig3}(a) and (b) are introduced. Although 
the rhombuses in Figure~\ref{fig3} have respectively the same shape as the rhombuses 
in Figure~\ref{fig1}, their matching conditions are different. For the two rhombuses in 
Figure~\ref{fig3}, the vertices are not colored, and all edges have orientations. To 
obtain the Penrose tiling (Figure~\ref{fig3}(c)) that is generated by the rhombuses 
in Figure~\ref{fig3}(a) and (b), the lengths and orientations of the edges of the 
rhombuses must match. On the other hand, the rhombus in Figure~\ref{fig3}(b) can 
independently generate a periodic tiling, as shown in Figure~\ref{fig3}(d), according 
to a matching condition. That is, although the two prototiles in Figure~\ref{fig3}(a) 
and (b) can generate a nonperiodic tiling, one of them can also generate a periodic 
tiling. Therefore, the set of the rhombuses in Figure~\ref{fig3}(a) and (b) is not an 
aperiodic set of prototiles. Based on the above discussion, we note that a set 
of prototiles of the Penrose tiling is not necessarily an aperiodic set of prototiles.

\renewcommand{\figurename}{{\small Figure.}}

\subsection{Aperiodic Tiling Due to Ammann That Uses Three 
Unmarked Convex Polygons as Prototiles}
\label{subsection2_2}

\begin{figure}[htbp]
  \centering\includegraphics[width=12.8cm,clip]{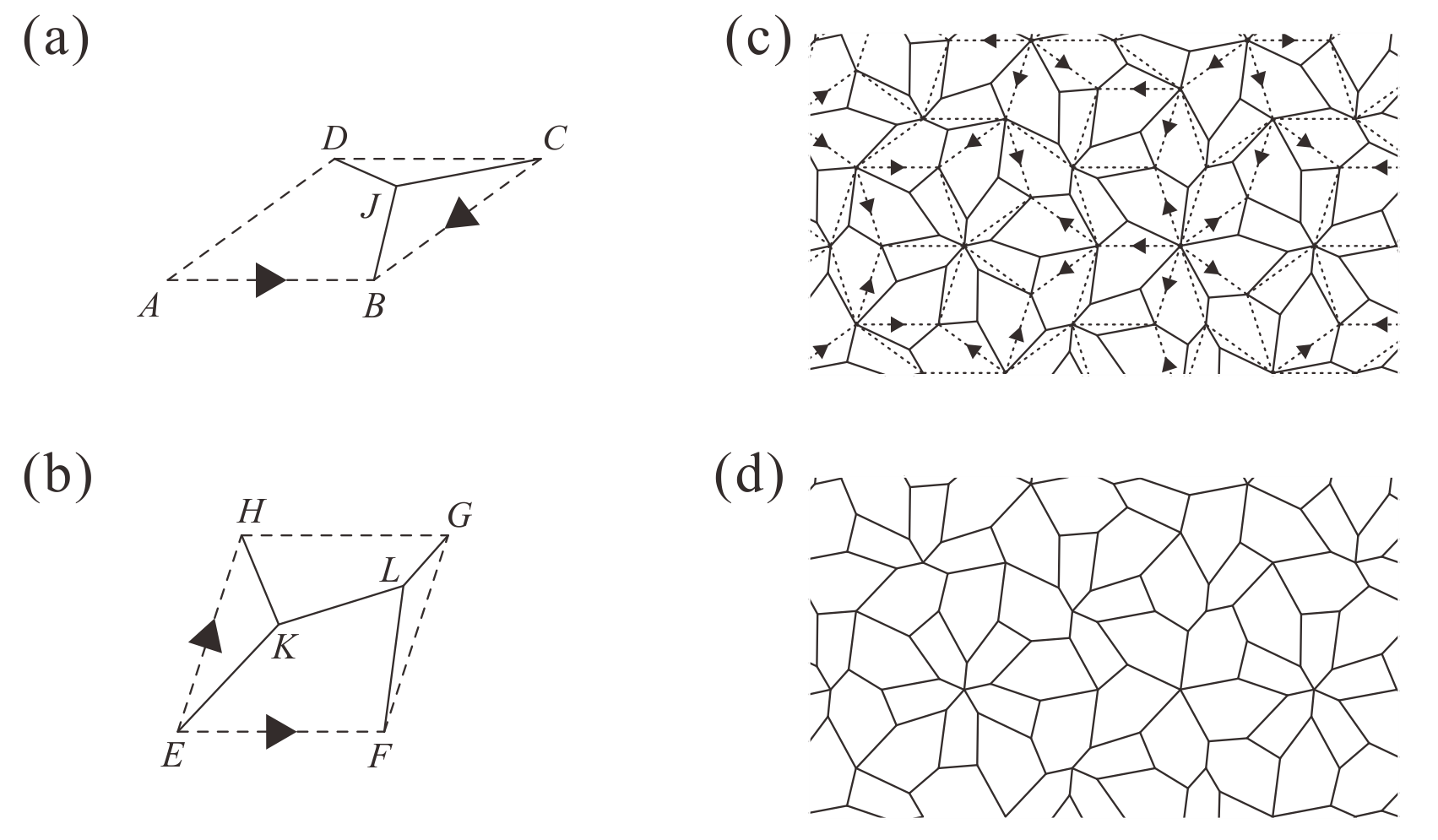}
  \caption{{\small 
Aperiodic tiling due to Ammann, which uses three unmarked 
convex polygons as prototiles---one hexagon and two pentagons. 
This tiling is a recomposition by Penrose rhombuses and is obtained 
from it using the markings shown in (a) and (b). These markings are 
completely determined by the choice of point $J$, as we must have 
\textit{GL}~=~\textit{DJ}, \textit{FL}~=~\textit{CJ}~=~\textit{EK}, 
and \textit{HK}~=~\textit{BJ}. 
\cite{G_and_S_1987} points out that the aperiodicity $J$ must 
be chosen so that \textit{DJ}, \textit{CJ}, \textit{BJ}, and \textit{KL} 
are of different lengths. \label{fig4}
}}
\end{figure}

In this subsection, we introduce the prototiles and tilings relevant to 
Problem (ii) on aperiodic tiling mentioned in Section~\ref{section1}. Problem (ii) 
is the problem of finding the smallest sets of aperiodic prototiles, each of which 
is an unmarked convex polygon (i.e., a convex polygon with no matching 
condition). Ammann has produced a remarkable example of such a set by 
recomposition from the set of rhombuses in Figure~\ref{fig1}(a) and 
(b)~\cite{G_and_S_1987}. If the inside of each rhombus in Figure ~\ref{fig1}(a) 
and (b) is divided as shown in Figure ~\ref{fig4}(a) and (b), the Penrose 
tiling in Figure~\ref{fig1}(c) becomes the tiling shown in Figure~\ref{fig4}(c). 
If the edges (broken lines) of the rhombuses are eliminated from the tiling 
in Figure~\ref{fig4}(c), the aperiodic tiling with one unmarked convex hexagon 
and two convex pentagons is obtained (see Figure~\ref{fig4}(d)). 
That is, the set of one convex hexagon and two convex pentagons in the 
edge-to-edge tiling in Figure~\ref{fig4}(d) is the set of three convex polygons 
that are aperiodic with no matching condition on the edges.

\renewcommand{\figurename}{{\small Figure.}}

\section{Proof of Theorem~\ref{t2}}
\label{section3}

A unit that can only generate a periodic tiling by translation is called a 
\textit{fundamental region}.

\begin{figure}[htbp]
 \centering\includegraphics[width=12.8cm,clip]{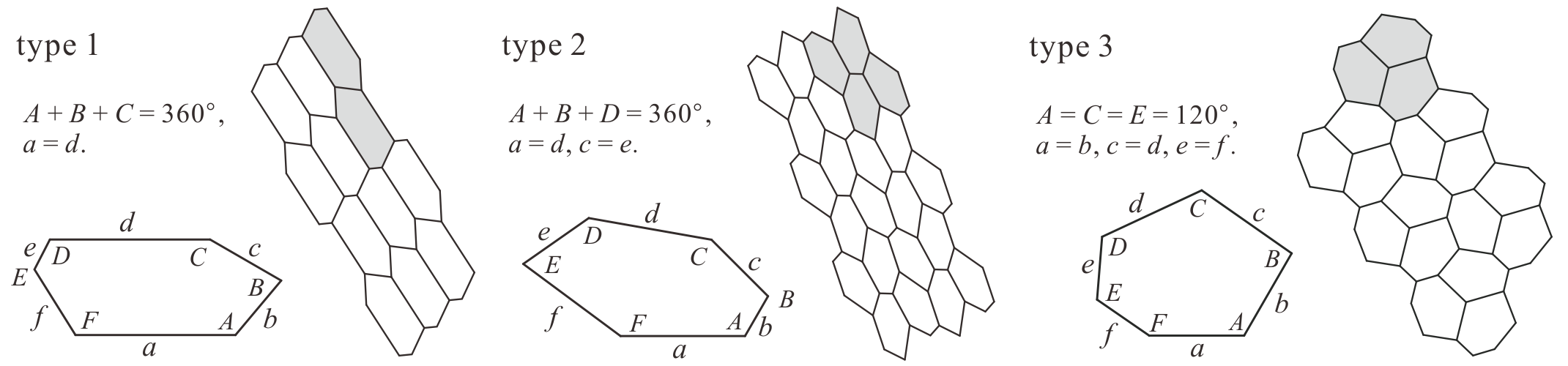}
  \caption{{\small 
Three types of convex hexagonal tiles. If a convex hexagon can 
generate a monohedral tiling, it belongs to at least one of types 1--3. 
The pale gray hexagons in each tiling indicate a fundamental region.}
\label{fig5}
}
\end{figure}

No convex polygon with seven or more edges can generate a monohedral 
tiling~\cite{Gardner_1975a, G_and_S_1987, Hallard_1991, Kershner_1968, Klamkin_1980, 
Reinhardt_1918}. All convex hexagons that can generate a monohedral tiling are 
categorized into three types\footnote{The classification of types of convex 
polygonal tiles is based on the essentially different properties of polygons. The 
classification problem of types of convex polygonal tiles and the classification 
problem of polygonal tilings are quite different.} and can generate a periodic 
edge-to-edge tiling (see Figure~\ref{fig5})~\cite{Gardner_1975a, G_and_S_1987, 
Hallard_1991, Kershner_1968, Reinhardt_1918}. 
Note that, a regular convex hexagon belongs to all of these three types.

\renewcommand{\figurename}{{\small Figure.}}

\begin{figure}[!p]
  \centering\includegraphics[width=12.8cm,clip]{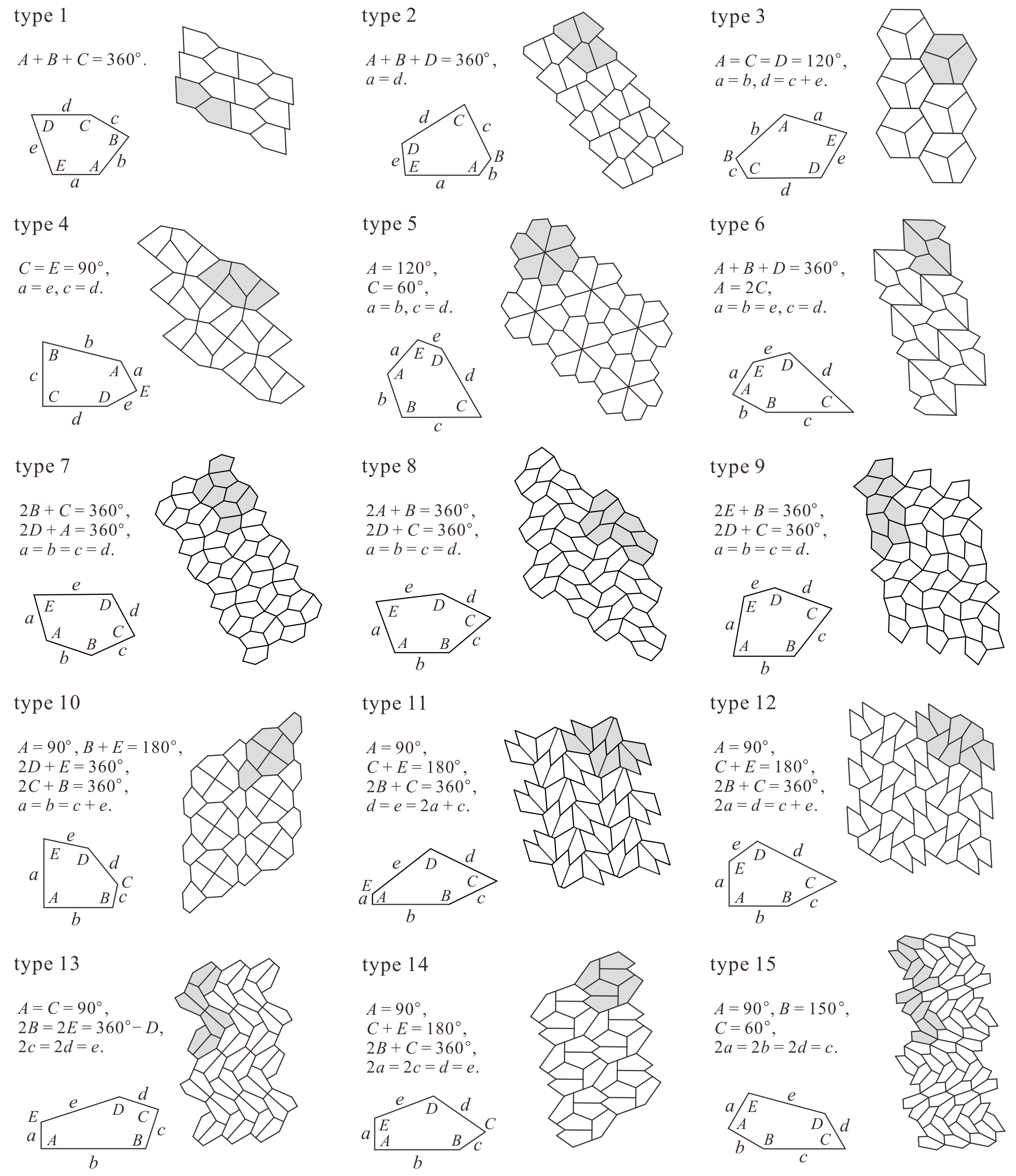}
  \caption{{\small 
Fifteen types of convex pentagonal tiles.
Each of the convex pentagonal tiles is defined by some conditions between 
the lengths of the edges and the magnitudes of the angles, but some degrees 
of freedom remain. For example, a convex pentagonal tile belonging to type 1 
satisfies that the sum of three consecutive angles is equal to $360^ \circ$. 
This condition for type 1 is expressed as $A + B + C = 360^ \circ$ in this figure. 
The pentagonal tiles of types 14 and 15 have one degree of freedom, that of 
size. For example, the value of $C$ of the pentagonal tile of type 14 is 
$\cos ^{ - 1}((3\sqrt {57} - 17) / 16) \approx 1.2099\;\mbox{rad} \approx 
69.32^ \circ$.}
\label{fig6}
}
\end{figure}

\begin{figure}[htbp]
  \centering\includegraphics[width=12.8cm,clip]{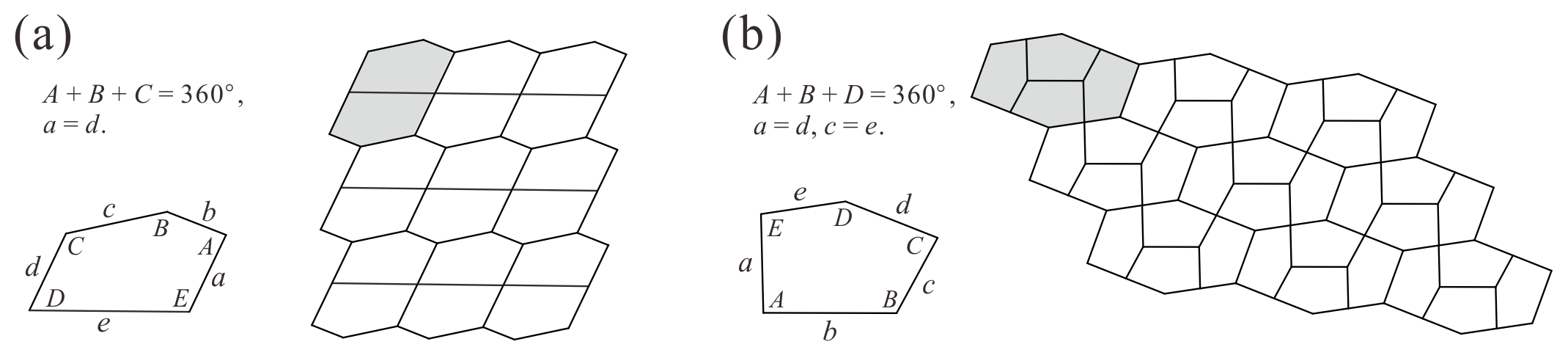}
  \caption{{\small 
Examples of edge-to-edge tilings by convex pentagonal tiles that 
belong to type 1 or type 2. The pale gray pentagons in each tiling indicate 
the fundamental region. (a) Convex pentagonal tiles that belong to type 1. 
(b) Convex pentagonal tiles that belong to type 2~\cite{Sugimoto_2012a, 
Sugimoto_2015}.}
\label{fig7}
}
\end{figure}

At present, essentially 15 different types of convex pentagonal tiles are 
known (see Figure~\ref{fig6}), but it is not known whether this list is 
complete~\cite{Gardner_1975a, G_and_S_1987, Hallard_1991, Kershner_1968, 
Mann_2015, Reinhardt_1918, Schatt_1978, 
Stein_1985, Sugimoto_2012a, Sugimoto_2012b, Sugimoto_2015, Sugimoto_III, 
Sugimoto_NoteTP, Sugi_Ogawa_2006, Wells_1991, wiki_pentagon_tiling}. 
There are edge-to-edge and non-edge-to-edge tilings by convex 
pentagonal tiles of the 15 types. A convex pentagonal tile belonging only 
to types 3, 10, 11, 12, 13, 14, or 15 in Figure~\ref{fig6} cannot generate 
an edge-to-edge tiling. In Figure~\ref{fig6}, non-edge-to-edge tilings 
are shown for type 1 or type 2 as general representative tilings, but 
the set of convex pentagonal tiles of type 1 or type 2 contains convex 
pentagonal tiles that can generate edge-to-edge tilings (see 
Figure~\ref{fig7}). The remaining six types can generate an edge-to-edge 
tiling~\cite{Sugimoto_2012a, Sugimoto_2015}. Because all convex pentagons 
that can generate an edge-to-edge monohedral tiling belong to at least one 
of eight types (i.e., types 1, 2, and 4--9 in Figure~\ref{fig6}) among the 15 
types~\cite{Bagina_2011, Sugimoto_2012b, Sugimoto_III, Sugimoto_NoteTP}, 
they must be able to generate a periodic tiling.

\renewcommand{\figurename}{{\small Figure.}}

\renewcommand{\figurename}{{\small Figure.}}

Any convex quadrilateral can generate a periodic edge-to-edge tiling by 
using a fundamental region of a convex hexagon with edges that have parallel 
opposite sides and are equal in length, which is made by the two 
quadrilaterals. Any triangle can generate a periodic edge-to-edge tiling by 
using a fundamental region of a parallelogram that is made by the two 
triangles. Thus, we obtain Theorem~\ref{t2}.
\hspace{3cm} $\square$

\section{Remarks}
\label{section4}

From the known results of convex polygonal tiles, if there is a single 
convex polygon such as an aperiodic prototile with no matching condition on 
the edges, the tile is a convex pentagonal tile that can generate an 
aperiodic non-edge-to-edge tiling. Further, if the complete list of all 
types of convex pentagonal tiles is obtained and if all convex pentagonal 
tiles in the complete list can generate a periodic tiling, no single convex 
polygon can be an aperiodic prototile.

\bigskip
\noindent
{\textbf{Acknowledgments.}
The authors would like to give heartful thanks to Emeritus 
Professor Hiroshi Maehara, University of the Ryukyus, and 
Professor Shigeki Akiyama, University of Tsukuba, whose 
enormous support and insightful comments were invaluable.

\end{document}